\documentclass[smallcondensed,final,natbib]{svjour3}
\usepackage{graphicx}
\usepackage{amsfonts,amsmath,amssymb}
\usepackage{mathptmx}
\usepackage{bm}
\usepackage[colorlinks,linkcolor=blue,citecolor=blue,urlcolor=blue]{hyperref}
\usepackage[ruled]{algorithm}
\usepackage[noend]{algpseudocode}
\alglanguage{pseudocode}

\newcommand{\x}{\vec{x}}
\newcommand{\z}{\vec{z}}
\renewcommand{\u}{\vec{u}} 
\newcommand{\vlambda}{\bm{\lambda}} 
\newcommand{\tL}{\tens{L}}
\newcommand{\prox}{\mathrm{prox}}
\providecommand{\norm}[1]{\lVert#1\rVert}

\title{Block-Simultaneous Direction Method of Multipliers}
\subtitle{A proximal primal-dual splitting algorithm for nonconvex problems with multiple constraints}
\author{Fred Moolekamp \and Peter Melchior}
\institute{
Fred Moolekamp
\at
Department of Astrophysical Sciences, Princeton University, Princeton, NJ 08544, USA\\
\email{fredem@princeton.edu}
\and
Peter Melchior
\at
Department of Astrophysical Sciences, Princeton University, Princeton, NJ 08544, USA\\
\email{peter.melchior@princeton.edu}
}

\begin{document}
\maketitle

\begin{abstract}
We introduce a generalization of the linearized Alternating Direction Method of Multipliers to optimize a real-valued function $f$ of multiple arguments with potentially multiple constraints $g_\circ$ on each of them.
The function $f$ may be nonconvex as long as it is convex in every argument, while the constraints $g_\circ$ need to be convex but not smooth.
If $f$ is smooth, the proposed Block-Simultaneous Direction Method of Multipliers (bSDMM) can be interpreted as a proximal analog to inexact coordinate descent methods under constraints.
Unlike alternative approaches for joint solvers of multiple-constraint problems, we do not require linear operators $\tL$ of a constraint function $g(\tL\ \cdot)$ to be invertible or linked between each other.
bSDMM is well-suited for a range of optimization problems, in particular for data analysis, where $f$ is the likelihood function of a model and $\tL$ could be a transformation matrix describing e.g. finite differences or basis transforms.
We apply bSDMM to the Non-negative Matrix Factorization task of a hyperspectral unmixing problem and demonstrate convergence and effectiveness of multiple constraints on both matrix factors.
The algorithms are implemented in {\sf python} and released as an open-source package.
\end{abstract}

\section{Introduction}

In this paper we seek to numerically
\begin{equation}
\label{eq:bsdmm}
\text{minimize}\ 
f\left(\x_{1},\dots,\x_{N}\right) + \sum_{j=1}^N \sum_{i=1}^{M_j}g_{ij}\left(\tL_{ij}\x_{j}\right),
\end{equation}
where $f:\ \mathbb{R}^{d_1}\times\dots\times\mathbb{R}^{d_N}\rightarrow \mathbb{R}$ is a potentially nonconvex function that is a closed proper convex function in each of $N$ independent arguments, 
and $g_{ij}: \mathbb{R}^{d_{ij}}\rightarrow \mathbb{R}$ are convex functions that encode $M_j$ constraints for each of those variables after they are mapped by linear operators $\tL_{ij}$.%
\footnote{Throughout this work, indices denote different variables or constraints, not elements of vectors or tensors.}

If $f$ and $g_\circ$ were all smooth, we could directly apply conventional gradient methods to solve the problem, but in many cases this is not possible. Examples are projections onto the positive orthant of $\mathbb{R}^{d_j}$ or regularization with the $\ell_1$ norm.
Instead we will access the functions through their proximal operators, defined as
\begin{equation}
\label{eq:proximal}
\prox_{\lambda f}\left(\vec{v}\right)\equiv\underset{\x}{\textrm{argmin}}\left\lbrace f\left(\x\right)+\frac{1}{2\lambda}\norm{\x-\vec{v}}_{2}^{2}\right\rbrace
\end{equation}
with a scaling parameter $\lambda$.
Their primary purpose is to turn any convex function $f$ into a strongly convex function even if $f$ may assume infinite values.
The minimization in \autoref{eq:proximal} may seem complicated, but in many cases exact and simple solvers exists.
For instance, if $f(\x)$ is the indicator function of the closed convex set $\mathcal{C}$, $\prox_f$ is simply the Euclidean projection operator onto $\mathcal{C}$.
Whether \autoref{eq:bsdmm} can efficiently be solved thus depends on the cost of evaluating the proximal operators involved. 
For more details and various interpretations of proximal operators we refer to \citet{Parikh2014} and references therein.

Several proximal algorithms exist to minimize a function $f$ of a single argument. 
If $f$ is smooth, the proximal gradient method provides minimization of $f(\x) + g(\x)$ with the \emph{forward-backward} scheme, where at iteration $k$ a step in the direction of $\nabla f$ is followed by the application of $\prox_g$:
\begin{equation}
\label{eq:pgm}
\x^{k+1} := \prox_{\lambda^k g} \left(\x^k - \lambda^k \nabla f(\x^k)\right).
\end{equation}
If the step size is $\lambda \in (0, 1/L]$ with $L$ being the Lipschitz constant of $\nabla f$, the convergence rate is $O(1/k)$, which can further be accelerated \citep{Nesterov2013}. 
We are particularly interested in introducing linear operators, which encode typical image processing operations, e.g. finite differences,  smoothing, or basis transforms.
A well-known approach to such a situation is the Alternating Direction Method of Multipliers (ADMM), which we review in \autoref{sec:ADMM}; several such constraints can be imposed with the Simultaneous Direction Method of Multipliers (SDMM, \autoref{sec:sdmm}).
We will introduce previous works in the relevant sections and state here only our main contributions:
We 1) generalize SDMM to cases where none of the $\tL_{ij}$ need to have full rank or be linked between different constraints, and 2) extend SDMM to nonconvex functions $f$ of several arguments $\x_j$, resulting in a proximal variant of inexact block optimization methods (\autoref{sec:bsdmm}).
As an example of the proposed algorithm, in \autoref{sec:nmf} we implement a solver for the Non-negative Matrix Factorization problem, where we allow for an arbitrary number of constraints on either of the matrix factors. We demonstrate convergence of the primal and dual variables as well as the effectiveness of multiple constraints in the example case.
We conclude in \autoref{sec:conclusions}.

\section{Generalizing ADMM}

\subsection{Linearized Alternating Direction Method of Multipliers}
\label{sec:ADMM}

We start by introducing the well-known Alternating Direction Method of Multipliers \citep[ADMM,][]{Gabay1976, Glowinski1975, Eckstein1992} in the notation that will be used throughout the paper.
It is applicable when $N=1$, i.e. to
\begin{equation}
\label{eq:ADMM}
\underset{\x_{1}}{\text{{minimize}}\:}f\left(\x_1\right)+g_1\left(\tL_{11}\x_1\right).
\end{equation}
We first re-write \autoref{eq:ADMM} in consensus form as 
\begin{equation}
\begin{array}{cc}
\textrm{{minimize}} & f\left(\x_{1}\right)+g_{1}\left(\z_{11}\right)\\
\textrm{{subject\,to}} & \tL_{11}\x_{1}-\z_{11}=0.
\end{array}
\end{equation}
The central idea is to split the optimization in two separate tasks: one that minimizes $f$ and another than satisfies $g_1$ by introducing the auxiliary variable $\z_{11}$ \citep{Douglas1956}.
This can be done by introducing Lagrange multipliers for each constraint, plus a quadratic term that is of critical importance for the convergence of the algorithm when $g_1$ is not strongly convex. 
The resulting Augmented Lagrangian for \autoref{eq:ADMM} is 
\begin{equation}
\mathcal{{L}}\left(\x_{1},\z_{11},\vlambda_{11}\right)=f\left(\x_{1}\right)+g_{1}\left(\z_{11}\right)+\vlambda_{11}^{\top}\left(\tL_{11}\x_{1}-\z_{11}\right)+\frac{1}{2\rho_{11}}\norm{\tL_{11}\x_{1}-\z_{11}}_{2}^{2},
\end{equation}
where $\rho_{11}\in\mathbb{{R}}>0$ and $\vlambda_{11}\in\mathbb{R}^{d_1}$. 
We then need to find unique minimizers of $\mathcal{{L}}$ with respect its variables.
In an iterative sequence one would formally need to update both $\x_1$ and $\z_{11}$ simultaneously, which is often numerically difficult. 
\citet{Chen1994} demonstrated that it is sufficient to update $\x_1$ first, then use the updated value of $\x_1$ for the $\z_{11}$-update, and then both values for the last update of $\vlambda_{11}$:
\begin{equation}
\begin{array}{ccl}
\x_{1}^{k+1} & := & \underset{\x_{1}}{\text{argmin}}\left\lbrace f\left(\x_{1}\right)+\vlambda_{11}^{k\top}\tL_{11}\x_{1}+\frac{1}{2\rho_{11}}\norm{\tL_{11}\x_{1}-\z_{11}^{k}}_{2}^{2}\,\right\rbrace\\
\z_{11}^{k+1} & := & \underset{\z_{11}}{\text{{argmin}}}\left\lbrace g_{1}\left(\z_{11}\right)-\vlambda_{11}^{k\top}\z_{11}+\frac{1}{2\rho_{11}}\norm{\tL_{11}\x_{1}^{k+1}-\z_{11}}_{2}^{2}\,\right\rbrace\\
\vlambda_{11}^{k+1} & := & \vlambda_{11}^{k}+\frac{1}{\rho_{11}}\left(\tL_{11}\x_{1}^{k+1}-\z_{11}^{k+1}\right).
\end{array}
\end{equation}
We can simplify the $\z_{11}$-update because we can add or subtract any terms independent of $\z_{11}$:
\begin{equation}
\label{eq:admm_zu}
\begin{array}{ccl}
\z^{k+1} & := & \underset{\z_{11}}{\text{{argmin}}}\left\lbrace g_{1}\left(\z_{11}\right)-\frac{1}{2\rho_{11}}\norm{\tL_{11}\x_{1}^{k+1}+\u_{11}^{k}-\z_{11}}_{2}^{2}\,\right\rbrace = \prox_{\rho_{11} g_{1}}\left(\tL_{11}\x_{1}^{k+1}+\u_{11}^{k}\right)\\
\u_{11}^{k+1} & := & \u_{11}^{k}+\tL_{11}\x_{1}^{k+1}-\z_{11}^{k+1}
\end{array}
\end{equation}
where we have introduced a scaled variable $\u_{11}^{k}\equiv\rho_{11}\vlambda_{11}^{k}$
and applied \autoref{eq:proximal}.

One cannot solve for $\x_{1}^{k+1}$ in the same way because of the presence of the linear operator $\tens{L}_{11}$.
\citet{Stephanopoulos1975} showed that linearizing $\frac{1}{2\rho_{11}}\norm{\tL_{11}\x_{1}-\z_{11}^{k}}_{2}^{2}$
at the current-iteration $\x^{k}_1$ is a practical solution that preserves the general convergence of the algorithm while providing a separable update sequence: 
\begin{equation}
\begin{array}{ccl}
\x_{1}^{k+1} & := & \underset{\x_{1}}{\text{argmin}}\left\lbrace f\left(\x_{1}\right)+\vlambda_{11}^{k\top}\tL_{11}\x_{1}+\frac{1}{\rho_{11}} \tL_{11}^{\top}\left(\tL_{11}\x_{1}^{k}-\z_{11}^{k}\right)\x_{1}+\frac{1}{2\mu_1}\norm{\x_{1}-\x_{1}^{k}}_{2}^{2}\,\right\rbrace\\
 & = & \prox_{\mu_1 f}\left(\x_{1}^{k}-\frac{\mu_1}{\rho_{11}} \tL_{11}^{\top}\left(\tL_{11}\x_{1}^{k}-\z_{11}^{k}+\u_{1}^{k}\right)\right),
\end{array}
\end{equation}
where we have introduced the parameter $\mu_1$ with $0<\mu_1\leq\rho_{11}/||\tL_{11}||_\mathrm{s}^2$.%
\footnote{We use $\norm{\cdot}_\mathrm{s}$ to denote the spectral norm, $\norm{\cdot}_2$ for the element-wise $\ell_2$ norm of vectors and tensors.}
The algorithm for a single variable $\x_1$ and a single constraint $g_1(\tens{L}_{11}\x_1)$ is also known as \emph{split inexact Uzawa method} \citep[e.g.][]{Esser2010, Parikh2014} and listed as \autoref{alg:admm}.

\begin{algorithm}[t]
\caption{Linearized ADMM}
\label{alg:admm}
\begin{algorithmic}[1]
\Procedure{ADMM}{$\x_1,\ \mu_1,\ \rho_{11},\ \tL_{11}$}
\State $\x_1^1 \gets \x_1;\ \z_{11}^1 \gets \tL_{11}\x_1;\ \u_{11}^1 \gets \vec{0}$
\For{$k = 1,2,\dots$}
    \State $\x_1^{k+1} \gets \prox_{\mu_1 f}\left(\x_{1}^{k}-\mu_1/\rho_{11} \tL_{11}^{\top}\left(\tL_{11}\x_{1}^{k}-\z_{11}^{k}+\u_{1}^{k}\right)\right)$
    \State $\z_{11}^{k+1} \gets \prox_{\rho_{11} g_{1}}\left(\tL_{11}\x_{1}^{k+1}+\u_{11}^{k}\right)$
    \State $\u_{11}^{k+1} \gets \u_{11}^{k}+\tL_{11}\x_{1}^{k+1}-\z_{11}^{k+1}$
    \If{$\norm{\vec{r}_{11}^{k+1}}_2 \leq \epsilon^\mathrm{pri} \wedge \norm{\vec{s}_{11}^{k+1}}_2 \leq \epsilon^\mathrm{dual}$}
     break
    \EndIf
\EndFor
\EndProcedure
\end{algorithmic}
\end{algorithm}

Following \citet{Boyd2011}, we implement stopping criteria based on primal residual $\vec{r}_{11}^{k+1}=\tens{L}_{11} \x_1^{k+1}-\z_{11}^{k+1}$ and the dual residual $\vec{s}_{11}^{k+1}=\frac{1}{\rho_{11}} \tens{L}_{11}^\top\left(\z_{11}^{k+1}-\z_{11}^{k} \right)$.
To assess primal and dual feasibility, we require
\begin{equation}
\label{eq:residuals}
\begin{split}
\norm{\vec{r}_{11}^{k+1}}_2 \leq \epsilon^\mathrm{pri} &\equiv \sqrt{p}\, \epsilon^\mathrm{abs} + \epsilon^\mathrm{rel} \max\lbrace \tens{L}_{11} \x_1^{k+1}, \z_{11}^{k+1}\rbrace\ \ \text{and}\\
\norm{\vec{s}_{11}^{k+1}}_2 \leq \epsilon^\mathrm{dual} &\equiv \sqrt{n}\, \epsilon^\mathrm{abs} + \epsilon^\mathrm{rel} / \rho_{11}\, \tens{L}_{11}^\top \u_{11}^{k+1},
\end{split}
\end{equation}
where $p$ and $n$ are the number of elements in $\z_{11}$ and $\x_1$, respectively.
The error thresholds $\epsilon^\mathrm{abs}$ and $\epsilon^\mathrm{rel}$ can be set at suitable values, depending on the precision and runtime constraints of the application.

\subsection{Several constraint functions}
\label{sec:sdmm}

It is often necessary to impose several constraints simultaneously.
\cite{Condat2013} proposed primal-dual split algorithms that can solve a restricted version of \autoref{eq:bsdmm}, in which $f(\x)$ is convex and differentiable with a Lipschitz-continuous gradient, and only two additional, potentially non-smooth functions are present,  $g(\x)$ and $h(\tL \x)$.
More generally applicable, \cite{Combettes2009} introduced the Simultaneous Direction Method of Multipliers (SDMM) to
\begin{equation}
\label{eq:sdmm_orig}
\underset{\x_{1}}{\text{{minimize}}}\sum_{i=1}^{M_1}g_{i}\left(\tL_{i1}\x_{1}\right),
\end{equation}
for which $g_i$ only have to be convex.
This corresponds to $N=1$ and $f\left(\x_{1}\right)=g_1(\x_1)$ with $L_{11}=\tens{I}$ in our notation, and enables the adoption of an arbitrary number of constraint functions, which is expressly what we seek.
However, their algorithm requires that $\tens{Q}=\sum_{i=1}^{M_1}\tL_{i1}^{\top}\tL_{i1}$ be invertible, a limitation that is too restrictive for many problems of interest.
We can dispense with this requirement by adopting the same linearization strategy as before with the ADMM.
Such a strategy is sensible if at least one function $g_l$ in \autoref{eq:sdmm_orig} can act as $f$ in \autoref{eq:bsdmm}.%
\footnote{While it is always possible to reformulate the problem thusly because we can set $f(\x_1) = g_l(\tL_{j 1} \x_1)$ for any $l$, it may render inefficient the minimization of $f$ by means of an proximal operator.
This is the limitation of the algorithm we derive in this section.}

Without loss of generality, we take $l=M_1$ and redefine $M_1 \rightarrow M_1 - 1$ in \autoref{eq:sdmm_orig} to maintain our notation.
We introduce $M_1$ primal variables $\z_{i1}$ and solve the problem in consensus form:
\begin{equation}
\label{eq:sdmm}
\begin{array}{cl}
\underset{\x_{1}}{\textrm{minimize}} & f\left(\x_{1}\right)+\sum_{i=1}^{M_{1}}g_{i}\left(\z_{i1}\right)\\
\textrm{{subject\,to}} & \tL_{i1}\x_{1}-\z_{i1}=0\ \ \forall i\in\lbrace 1,\dots,M_1\rbrace.
\end{array}
\end{equation}
Updating the primal and dual variables  $\u_{i1}$ is exactly the same as in \autoref{eq:admm_zu}, 
\begin{equation}
\begin{array}{ccl}
\z_{i1}^{k+1} & := & \prox_{\rho_{i1} g_{i1}}\left(\tL_{i1}\x_{1}^{k+1}+\u_{i1}^{k}\right)\\
\u_{i1}^{k+1} & := & \u_{i1}^{k}+\tL_{i1}\x_{1}^{k+1}-\z_{i1}^{k+1},
\end{array}
\end{equation}
where each update can be performed in parallel, justifying the ``S'' in SDMM.
Because we isolated $f$ as a function of $\x_1$ alone, we can minimize $\mathcal{L}\left(\x_{1},\z_{11},\dots,\z_{M_{1}1},\vlambda_{11},\dots,\vlambda_{M_{1}1}\right)$ with respect to $\x_1$ with the update
\begin{equation}
\label{eq:sdmm_x}
\begin{array}{ccl}
\x_{1}^{k+1} & := & \underset{\x_{1}}{\text{{argmin}}}\left\lbrace f\left(\x_{1}\right)+\sum_{i=1}^{M_{1}}\left(\vlambda_{i1}^{k\top}\tL_{i1}\x_{1}+\frac{1}{2\rho_{i1}}\norm{\tL_{i1}\x_{1}-\z_{i1}^{k}}_{2}^{2}\right)\,\right\rbrace\\
 & = & \underset{\x_{1}}{\text{{argmin}}}\left\lbrace f\left(\x_{1}\right)+\sum_{i=1}^{M_{1}}\left(\vlambda_{i1}^{k\top}\tL_{i1}\x_{1}+\frac{1}{\rho_{i1}} \tL_{i1}^{\top}\left(\tL_{i1}\x_{1}^{k}-\z_{i1}^{k}\right)\x_{1}+\frac{1}{2\mu_1}\norm{\x_{1}-\x_{1}^{k}}_{2}^{2}\right)\,\right\rbrace\\
 & = & \prox_{\mu_1 f}\left(\x_{1}^{k}-\sum_{i=1}^{M_{1}}\mu_1/\rho_{i1}\tL_{i1}^{\top}\left(\tL_{i1}\x_{1}^{k}-\z_{i1}^{k}+\u_{i1}^{k}\right)\right),
\end{array}
\end{equation}
where we linearized in the second step and added a quadratic penalty to introduce $\prox_{\mu_1 f}$ in the third step. The parameters $\mu_1$ and $\rho_{i1}$ are bound by
\begin{equation}
\label{eq:sdmm_rho}
\rho_{i1}/\norm{\tL_{i1}}_\mathrm{s}^2 \geq \beta_\mathrm{S}\, \mu_1 \text{ with } 1 \leq \beta_\mathrm{S} \leq M_1.
\end{equation}
The parameter $\beta_\mathrm{S}$ is necessary to account for potentially correlated contributions of different  $g_i$ in \autoref{eq:sdmm_x}. 
If the $g_i$ are partially degenerate and one would adopt the na\"ive threshold $\beta_\mathrm{S}=1$, $\prox_{\mu_1 f}$ will lose its contracting property.
The most conservative option $\beta_\mathrm{S}=M_1$ will always lead to a convergent minimizer, even if $g_i=g\ \forall i$, albeit generally at the expense of reduced convergence speeds.

The linearized form of SDMM is listed in \autoref{alg:sdmm}.
We note that the resulting algorithm is not identical and thus not suited to the same set of problems as the original SDMM by \cite{Combettes2009} because of the isolation of $f(\x_1)$ in \autoref{eq:sdmm} that is not present in \autoref{eq:sdmm_orig}.
However, a strong similarity persists with the exception of Line 4 in \autoref{alg:sdmm}, we thus consider it appropriate to call this algorithm Linearized SDMM.

\begin{algorithm}[t]
\caption{Linearized SDMM}
\begin{algorithmic}[1]
\Procedure{SDMM}{$\x_1, \mu_1, [\rho_{11}\dots,\rho_{M_1 1}],\ [\tL_{11}\dots,\tL_{M_1 1}]$}
\State $\x_1^1 \gets \x_1;\ \z_{i1}^1 \gets \tL_{i1}\x_1\ \forall i \in \lbrace 1,\dots,M_1\rbrace;\ \u_{i1}^1 \gets \vec{0}\ \forall i \in \lbrace 1,\dots,M_1\rbrace$
\For{$k = 1,2,\dots$}
    \State $\x_1^{k+1} \gets\prox_{\mu_1 f}\left(\x_{1}^{k}-\sum_{i=1}^{M_{1}}\mu_1/\rho_{i1}\tL_{i1}^{\top}\left(\tL_{i1}\x_{1}^{k}-\z_{i1}^{k}+\u_{i1}^{k}\right)\right)$
    \For{$i = 1,\dots,M_1$}
    	\State $\z_{i1}^{k+1} \gets \prox_{\rho_{i1} g_{i1}}\left(\tL_{i1}\x_{1}^{k+1}+\u_{i1}^{k}\right)$
   	\State $\u_{i1}^{k+1} \gets \u_{i1}^{k}+\tL_{i1}\x_{1}^{k+1}-\z_{i1}^{k+1}$
   \EndFor
   \If{$\bigwedge_i\left\lbrace\norm{\vec{r}_{i1}^{k+1}}_2 \leq \epsilon^\mathrm{pri} \wedge \norm{\vec{s}_{i1}^{k+1}}_2 \leq \epsilon^\mathrm{dual}\right\rbrace$}
     break
    \EndIf
\EndFor
\EndProcedure
\end{algorithmic}
\label{alg:sdmm}
\end{algorithm}

\subsection{Nonconvex, multi-argument functions}
\label{sec:bsdmm}

We now seek to solve \autoref{eq:bsdmm}, which includes the treatment of the function $f$ having several arguments.
Even with our requirement that $f$ be a closed proper convex function in each of its arguments, $f$ itself is generally not convex
(we refer the reader to e.g. \citealt{Zhang2016} for cases when a constraint function $g$ is not convex).

The nonconvexity that arises in multi-argument functions has been addressed with an ADMM variant first by \citet{Hong2014}, who provide a provably convergent solution for
\begin{equation}
\label{eq:Hong14}
\begin{array}{cl}
\underset{\x_{1},\dots,\x_{N}}{\text{minimize}} & 
f\left(\x_{1},\dots,\x_{N}\right) + \sum_{j=1}^N g_{j}\left(\x_{j}\right)\\
\textrm{subject\,to} & \sum_j^N \tL_{j}\x_{j}=\vec{b},
\end{array}
\end{equation}
and later by \citet{Wang2015}, who solve
\begin{equation}
\begin{array}{cl}
\underset{\x_{1},\dots,\x_{N},\vec{y}}{\text{minimize}} & 
f\left(\x_{1},\dots,\x_{N},\vec{y}\right)\\
\textrm{subject\,to} & \sum_j^N \tL_{j}\x_{j} + \mathsf{B} \vec{y}=\vec{b}.
\end{array}
\end{equation}
While close to our problem in \autoref{eq:bsdmm}, especially the form of \autoref{eq:Hong14}, these approaches need the linking $\sum_j^N \tL_{j}\x_{j}=\vec{b}$ across the variables $\x_j$ (and $\vec{y}$ for \citealt{Wang2015}), which conflicts with our desire to have independent constraint functions.
In addition, inequality constraints like $\tL\x\leq0$ cannot be expressed in either of the forms above.

By now our strategy for solving \autoref{eq:bsdmm} with $N$ variables $\x_{j}$ should be apparent.
We first bring the problem into consensus form, i.e. we seek to
\begin{equation}
\begin{array}{cl}
\underset{\x_{1},\dots,\x_{N}}{\text{minimize}} & 
f\left(\x_{1},\dots,\x_{N}\right) + \sum_{j=1}^N \sum_{i=1}^{M_j}g_{ij}\left(\z_{ij}\right)\\
\textrm{subject\,to} & \tL_{ij}\x_{j}-\z_{ij}=0\ \ \forall j \in \lbrace 1,\dots,N\rbrace,\ i\in\lbrace 1,\dots,M_j \rbrace.
\end{array}
\end{equation}
This results in $\sum_{j=1}^{N}M_{j}$ primal and dual variables that are updated using 
\begin{equation}
\begin{array}{ccl}
\z_{ij}^{k+1} & := & \prox_{\rho_{ij} g_{ij}}\left(\tL_{ij}\x_{j}^{k+1}+\u_{ij}^{k}\right)\\
\u_{ij}^{k+1} & := & \u_{ij}^{k}+\tL_{ij}\x_{j}^{k+1}-\z_{ij}^{k+1}.
\end{array}
\end{equation}
Then we linearize the quadratic term $\frac{1}{2\rho_{ij}}\norm{\tL_{ij}\x_{j}-\z_{ij}^{k}}_{2}^{2}$, and utilize that $f$ is convex in every argument to solve for each $\x_j$ with a proximal operator as minimizer of $f$ wrt $\x_j$:
\begin{equation}
\begin{array}{ccl}
\x_{j}^{k+1} & := & \underset{\x_{j}}{\text{{argmin}}} \left\lbrace f\left(\x_{1},\dots,\x_{N}\right)+\sum_{i=1}^{M_{j}}\left(\vlambda_{ij}^{k\top}\tL_{ij}\x_{j}+\frac{1}{2\rho_{ij}}\norm{\tL_{ij}\x_{j}-\z_{ij}^{k}}_{2}^{2}\right)\,\right\rbrace\\
 & = & \underset{\x_{j}}{\text{{argmin}}} \left\lbrace f\left(\x_{1},\dots,\x_{N}\right)+\sum_{i=1}^{M_{1}}\left(\vlambda_{ij}^{k\top}\tL_{ij}\x_{1}+\frac{1}{\rho_{ij}} \tL_{ij}^{\top}\left(\tL_{ij}\x_{j}^{k}-\z_{ij}^{k}\right)\x_{j}+\frac{1}{2\mu_j}\norm{\x_{j}-\x_{j}^{k}}_{2}^{2}\right)\,\right\rbrace\\
 & = & \prox_{\mu_j f,j}\left(\x_{j}^{k}-\sum_{i=1}^{M_{j}}\frac{\mu_j}{\rho_{ij}}\tL_{ij}^{\top}\left(\tL_{ij}\x_{j}^{k}-\z_{ij}^{k}+\u_{ij}^{k}\right)\right),
\end{array}
\end{equation}
The entire algorithm, which we call bSDMM, is listed as \autoref{alg:bsdmm}.
If $f$ is separable, bSDMM amounts to $N$ independent, and thus parallelizable, solutions for each $\x_j$ using \autoref{alg:sdmm}.
If not, to the knowledge of the authors, general convergence guarantees do not exist since the solutions for each $\x_j$ may not be unique.
However, in the case of quadratic convex functions and of $N=2$ convergence can be guaranteed \citep{Grippo2000, Lin2007}.
If $f$ is smooth, block-wise optimization by means of successive proximal forward-backward steps is convergent even with Nesterov-type acceleration \citep{Razaviyayn2013, Xu2013}.
In that case, these optimization steps constitute a proximal variant of a block coordinate descent algorithm.
Our approach is different because it uses a primal-dual split instead of a direct constraint projection, which can only deal with one constraint per optimization variable, and the nonconvex ADMM variants of \citet{Hong2014} and \citet{Wang2015} because of the independence of the constraints.

Furthermore, because the updates do not necessarily yield the minimum of $f$ or the constraints in each step, e.g. for a single step of $\nabla_j f$, it is very efficient, but the limit point of the sequence is not necessarily a local minimum, just a stationary point.
As long as the approximations are sufficiently precise, ADMM is still convergent \citep{Eckstein1992, Eckstein2017}.
\citet{Berry2007} studied approximate solvers in several different applications and found acceptable results at a fraction of the computational cost.
We will inspect the convergence properties of the bSDMM algorithm with an suitable example in the next section.

As with the SDMM, the presence of several constraints $g_{ij}$ for a single function $f$ necessitates restraint when choosing $\rho_{ij}$ so as not to overwhelm $\prox_{\mu_j f,j}$:
\begin{equation}
\label{eq:bsdmm_rho}
\rho_{ij}/\norm{\tL_{ij}}_\mathrm{s}^2 \geq \beta_\mathrm{bS}\, \mu_j \text{ with } 1 \leq \beta_\mathrm{bS} \leq N\,M_j.
\end{equation}
However, because $f$ is now a function with several arguments, the parameter $\mu_j$ may change with every iteration.
For instance, if $f$ is smooth, $\mu_j$ is bound by the Lipschitz constant of $\nabla_j f\left(\x_{1},\dots,\x_{N}\right)$.
The bSDMM algorithm therefore requires the function $h\left(j;\,\x_{1},\dots,\x_{N}\right)$ to compute $\mu_j$, from which it will then determine $\rho_{ij}$ to satisfy \autoref{eq:bsdmm_rho} (lines 8 and 9 of \autoref{alg:bsdmm}).

\begin{algorithm}[t]
\caption{Block-SDMM}
\begin{algorithmic}[1]
\Procedure{bSDMM}{$[\x_1,\dots,\x_N],\ h,\ \beta_\mathrm{G},\ [\tL_{11}\dots,\tL_{M_N N}]$}
\For{$j=1,\dots,N$}
    \State $\x_j^1 \gets \x_j$
    \State $\z_{ij}^1 \gets \tL_{ij}\x_j\ \forall i \in \lbrace 1,\dots,M_j\rbrace$
    \State $\u_{ij}^1 \gets \vec{0}\ \forall i \in \lbrace 1,\dots,M_j\rbrace$
\EndFor
\For{$k = 1,2,\dots$}
    \For{$j=1,\dots,N$}
    	\State $\mu_j^{k+1} \gets h\left(j;\,\x_{1}^{k},\dots,\x_{N}^{k}\right)$
	\State $\rho_{ij}^{k+1} \gets \beta_\mathrm{G}\, \mu_j^{k+1} \norm{\tL_{ij}}_2^2$ 
        \State $\x_j^{k+1} \gets\prox_{\mu_j f,j}\left(\x_{j}^{k}-\sum_{i=1}^{M_{j}}\mu_j/\rho_{ij}\tL_{ij}^{\top}\left(\tL_{ij}\x_{j}^{k}-\z_{ij}^{k}+\u_{ij}^{k}\right)\right)$
        \For{$i = 1,\dots,M_j$}
    	    \State $\z_{ij}^{k+1} \gets \prox_{\rho_{i1} g_{i1}}\left(\tL_{i1}\x_{j}^{k+1}+\u_{ij}^{k}\right)$
   	    \State $\u_{ij}^{k+1} \gets \u_{ij}^{k}+\tL_{ij}\x_{j}^{k+1}-\z_{ij}^{k+1}$
	\EndFor
   \EndFor
      \If{$\bigwedge_{ij}\left\lbrace\norm{\vec{r}_{ij}^{k+1}}_2 \leq \epsilon^\mathrm{pri} \wedge \norm{\vec{s}_{ij}^{k+1}}_2 \leq \epsilon^\mathrm{dual}\right\rbrace$}
     break
    \EndIf
\EndFor
\EndProcedure
\end{algorithmic}
\label{alg:bsdmm}
\end{algorithm}

\section{Non-negative Matrix Factorization}
\label{sec:nmf}

An important application of bSDMM is Non-negative Matrix Factorization \citep[NMF,][]{Paatero1994}, which seeks to describe a data set $\tens{D}\in\mathbb{R}^{B \times L}$ of $L$-dimensional features that are observed $B$ times as a product of two non-negative matrices $\tens{A}\in\mathbb{R}^{B\times K}$ and $\tens{S}\in\mathbb{R}^{K\times L}$.
The idea is to reduce the dimensionality of the problem to $K$ prototypes, encoded in $\tens{S}$, whose sum generates the data in each observations with relative amplitudes encoded in $\tens{A}$.

We will adopt the Euclidean cost function, which corresponds to the negative log-likelihood under the assumption of standard Gaussian errors on each element of $\tens{D}$ (see \citet{Blanton2007,Zhu2016} for extensions to heteroscedastic errors).
The objective function is thus
\begin{equation}
f(\tens{A},\tens{S})=\norm{\tens{A}\cdot\tens{S} - \tens{D}}_2^2,
\end{equation}
and the non-negative constraints can be expressed as 
\begin{equation}
\label{eq:gplus}
g_+(\tens{A}) + g_+(\tens{S})\ \text{where}\ g_+(\tens{X}) = \begin{cases}
0 & \text{if}\ \tens{X}_{mn} \geq 0\ \forall m,n\\
\infty & \text{else}.
\end{cases}
\end{equation}
In the notation of \autoref{eq:bsdmm}, this corresponds to the minimally non-trivial case of $N=2$ and $M_j=1$, $\tL_{ij}=\tens{I}$ for $j={1,2}$.
The most basic NMF solver uses a ``multiplicative update'' (MU) rule that can be derived from a gradient descent argument \citep{Lee2001}.
However, MU has long been criticized for its often inferior convergence properties, which can be traced back to the implicit treatment of the  constraint.
Several alternative approaches have been brought forward to address the shortcomings of MU solvers, including the ability to impose constraints other than non-negativity \citep[e.g.][]{Berry2007}.
A comparison of these different approaches is not the focus of this work (see \citet{Xu2013} for a recent overview).
Instead, we demonstrate that bSDMM can successfully and efficiently impose several constraints on both matrix factors. 

To do so, we need the proximal-operator forms of the desired constraints.
In cases where the constraint is given by an indicator function of a convex set $\mathcal{C}$, the proximal operator is simply the projection operator onto $\mathcal{C}$ under the Euclidean norm, and the step size $\lambda$ is irrelevant.
For example, the non-negativity constraint becomes the (element-wise) projection onto the non-negative orthant:
\begin{equation}
\prox_{\lambda g_+} (\x) = \max(\bf{0},\x).
\end{equation}
Many other proximal operators can be evaluated analytically, e.g. for penalty functions involving $\ell_p$ norms, Total Variation, Maximum Entropy \citep{Combettes2009, Parikh2014}.
If $K$ is large, the data are noisy, or $B\ll L$, the NMF factors are generally degenerate.
Additional constraints then become necessary for reasonable results.

\subsection{Example: Hyperspectral unmixing}
Hyperspectral data are images takes of the same scene at several wavelengths ($B\sim100$), often beyond the range visible to humans.
This extended and more fine-grained spectral information can be used to robustly identify distinct components in the images.
If the spatial resolution of the imager is high enough, those components can be observed in their pure form, i.e. each pixel received contributions only from one component, but in general the spectral information of pixels is mixed.
Under the assumption that the components do not interfere, the so-called ``linear mixing model'', the NMF allows us to unmix the spectral contributions of each pixel, simultaneously inferring the pure spectrum of each of $K$ components, called ``endmembers'', as well as their amplitude in each pixel \citep{Berry2007, Jia2009, Gillis2014}.

\begin{figure}
    \includegraphics[width=\textwidth]{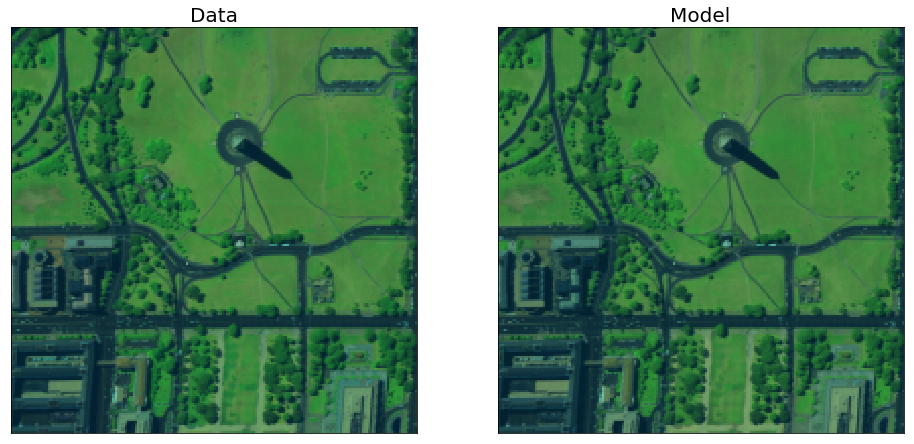}
    \caption{False-colored images of hyperspectral data and a four-component model comprising ``concrete'', ``soil'', ``vegetation'', and a spatially flat ``background''. The false-color image maps the sum of the first 50 wavelengths to blue, the next 50 to green, and the remaining 91 wavelengths to red.}
    \label{fig:false_color}
\end{figure}

\autoref{fig:false_color} shows a false-color image of the National Mall in Washington D.C. made from hyperspectral HYDICE \citep{Mitchell1995} data comprising $B=191$ wavelengths from 400 to 2475\,nm.%
\footnote{Data set obtained from \url{https://engineering.purdue.edu/~biehl/MultiSpec/}}
By our convention, the columns of matrix $\tens{A}$ describe the endmember spectra, and the rows of $\tens{S}$ the endmember amplitude per pixel.
Despite the large number of wavelengths, the problem is still strongly underconstrained even with a small number of endmembers.
We therefore impose several constraints:
\begin{itemize}
\item To prevent the degeneracy between $\tens{A}$ and $\tens{S}$ that stems from the transformation $(\tens{A},\tens{S})\rightarrow(\tens{A}\tens{Q},\tens{Q}^{-1}\tens{S})$ with an arbitrary invertible matrix $\tens{Q}$, we normalize the endmember spectra, i.e. the columns of $\tens{A}$. This normalization is different from the one usually adopted in hyperspectral unmixing applications, where the endmember amplitudes are normalized in each pixel, which results in endmembers being defined by both shape and amplitude of the spectrum. We prefer our approach because it maintains spectral similarity between regions of different brightness. 
\item The radiation recorded by the hyperspectral camera is a combination of light reflected off the ground and the atmosphere. Since we are interested in the former, and the latter is not expected to vary over the image, we add a ``background'' component that we constrain to be spatially flat. 
\item As the scene on the ground is mostly coherent over large areas, we  add a two-dimensional anisotropic total variation (TV) penalty \citep{Chambolle1997, Chambolle2004}. 
\end{itemize}
In summary, we minimize
\begin{equation}
\label{eq:unmixing}
\norm{\tens{A}\cdot\tens{S} - \tens{D}}_2^2 + g_+(\tens{A}) + g_+(\tens{S}) + g_\mathrm{norm}(\tens{L}_\mathrm{norm}\tens{A}) + g_\mathrm{bg}(\tens{S}) + \lambda \left(\norm{\tens{G}_x \tens{S}}_1 + \norm{\tens{G}_y \tens{S}}_1\right),
\end{equation}
where $g_+$ is given in \autoref{eq:gplus}, and $g_\mathrm{norm}$ and $g_\mathrm{bg}$ are additional indicator functions.
In detail, we combine the NMF fidelity term and the positivity constraints into forward-backward operators of the form of \autoref{eq:pgm}, one for $\tens{A}$ and one for $\tens{S}$.
This is the proximal-gradient technique for solving the NMF \citep[e.g.][]{Xu2013}, implemented in Line 10 of \autoref{alg:bsdmm}.

The normalization $\tens{L}_\mathrm{norm} \tens{A} = \tens{1}_K$ is achieved by $\tens{L}_\mathrm{norm}=\tens{1}_B^\top$, for which the proximal operator is simply a projection onto $\tens{1}_K$.

The background component requires that $\tens{S}_\mathrm{bg}=\mathrm{const}$ for all pixels. Since the proximal operator must yield the closest point on the submanifold in the Euclidean norm, it alters the background component row $\tens{S}_\mathrm{bg}\rightarrow \langle \tens{S}_\mathrm{bg} \rangle_L\, \tens{1}_L$, where the expectation value is carried out over all pixels, while leaving all other components unchanged.

For the TV penalty, we use the gradient operators in horizontal ($\tens{G}_x$) and vertical ($\tens{G}_y$) direction and make use of the analytic form of the proximal operator for the $l_1$ norm, the soft-thresholding operator \citep[e.g.][]{Combettes2005}. 
While we could adjust the penalty parameter $\lambda$ for the horizontal and vertical direction, as well as for every component, we have not found it necessary to explore that option.

The last four constraint and penalty functions are implemented in the SDMM fashion, giving rise to one auxiliary variable for $\tens{A}$ and three for $\tens{S}$. 

We run bSDMM with a TV penalty of $\lambda=10$ until feasibility with $\epsilon^\mathrm{rel}=0.01$ for primal and dual residuals is reached ($\epsilon^\mathrm{abs}$ was set to zero).
We chose the number of endmembers to be three, which we label as ``concrete'', ``soil'', and ``vegetation'', plus the flat background.%
\footnote{The choice of $K=4$ is somewhat arbitrary, and we have not attempted to find the optimal number of components since that is not the relevant aspect of this work.}

While random initialization of the matrices works reasonably well, we found better results when we initialize the spectra by a three-step approach: 
First, we determine the background spectrum by the minimum value over all pixels at a given wavelength.
Second, we select a reference pixel that, from its location in the scene, should to be a pure representation of one of the three other components, and subtract from that pixel the background spectrum.
Third, we normalize all spectra to sum up to one.
The initialization of $\tens{S}$ is less important.
We start with a zero matrix and then utilize bSDMM to make suitable updates given the initialized spectra.

\begin{figure}
    \includegraphics[width=\textwidth]{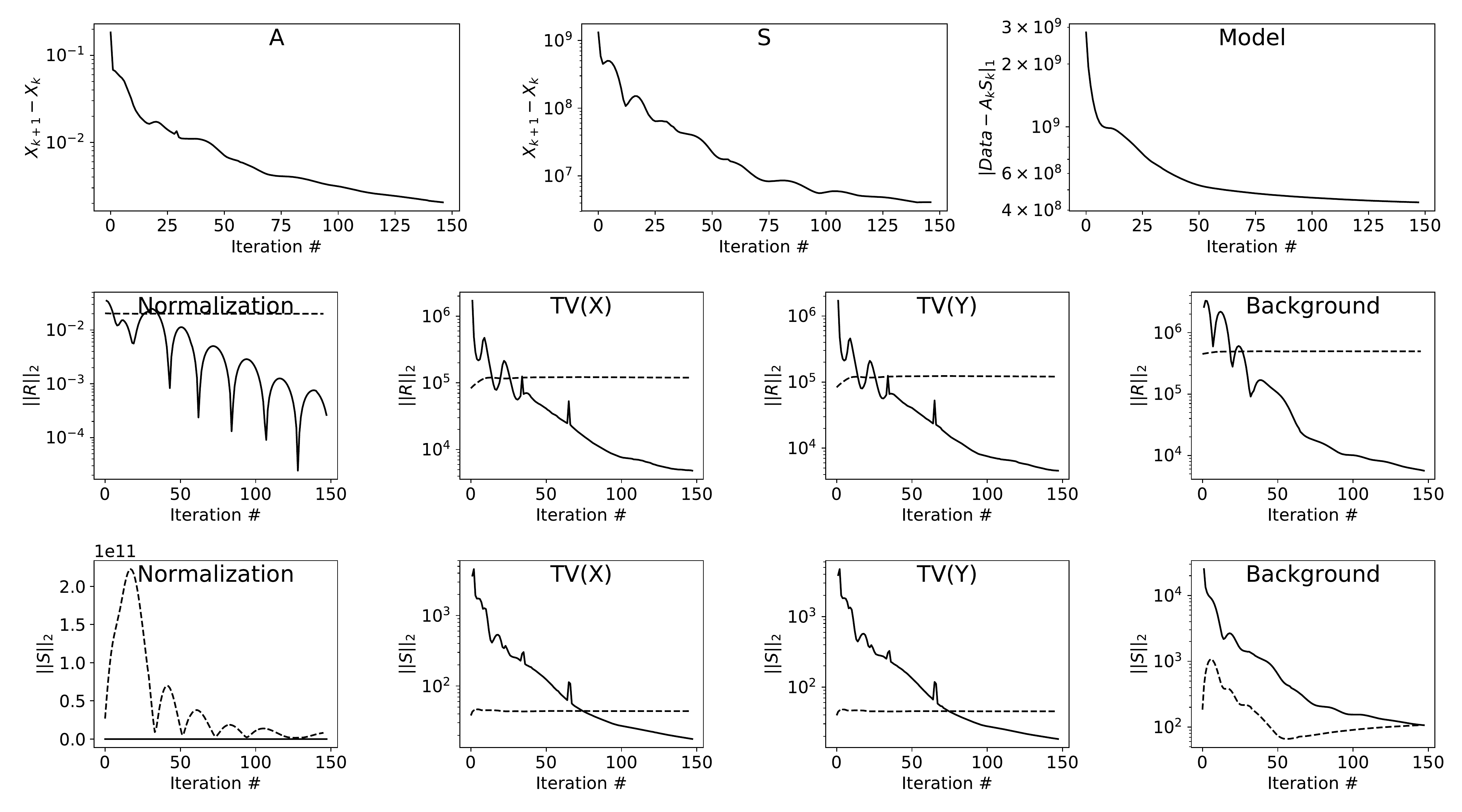}
    \caption{Convergence to a stationary point of the problem in \autoref{eq:unmixing}: $\tens{A}$, $\tens{S}$, and the model $\tens{A}\cdot\tens{S}$ (\emph{top});
    the primal residual $\vec{r}_{ij}^k$ (\emph{middle}) and the dual residual $\vec{s}_{ij}^k$ (\emph{bottom}) for the four constraint variables. 
    The limits $\epsilon^\mathrm{pri}$ and $\epsilon^\mathrm{dual}$ for primal and dual feasibility (cf. \autoref{eq:residuals}) are shown as dashed lines with $\epsilon^\mathrm{rel}=0.01$.}
    \label{fig:convergence}
\end{figure}

\begin{figure}
    \includegraphics[width=\textwidth]{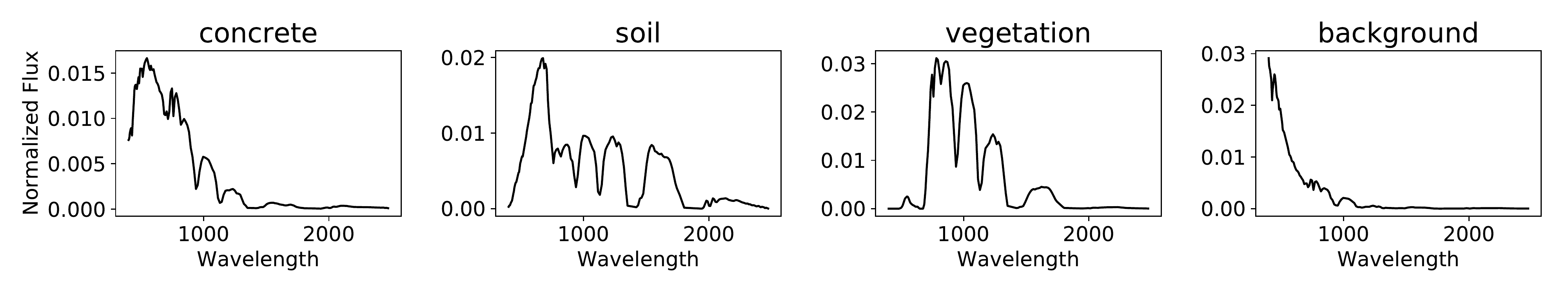}
    \caption{Endmember spectrum for every component. Each spectrum is normalized to unity. The labels are approximate descriptions given the regions in the image the endmembers mostly represent.}
    \label{fig:spec}
\end{figure}

\begin{figure}
    \includegraphics[width=\textwidth]{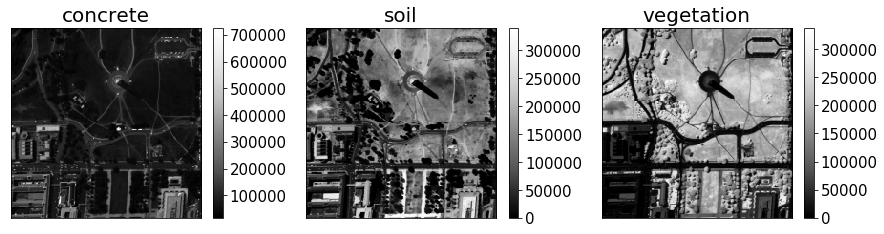}
    \caption{Intensity of the three spatially variable components ``concrete'', ``soil'', and ``vegetation''.
    Several features are prominent in the endmember intensities, such as trees, trails, and rooftops.
    These plots do not represent an endmember classification: because of our spectrum normalization, a region appears dark in these plots if it has a different spectrum than the endmember or if it reflects very little light (e.g. road surfaces).
    }
    \label{fig:intensity}
\end{figure}

The convergence of $\tens{A}$, $\tens{S}$, and the entire model $\tens{A}\cdot\tens{S}$ is shown in the top row of \autoref{fig:convergence},
primal and dual feasibility according to \autoref{eq:residuals} in the middle and bottom row.
The residual requirements are shown as dashled lines in \autoref{fig:convergence}, from which we can conclude that primal feasibility is generally achieved after about 30 iterations, while the dual feasibility of the background component requires almost 150 iterations.
The resulting endmember spectra and amplitudes are shown in \autoref{fig:spec} and \autoref{fig:intensity}.
The complete code to reproduce this hyperspectral unmixing example with the bSDMM-NMF approach is available at \url{https://github.com/fred3m/hyperspectral}.

\section{Conclusion}
\label{sec:conclusions}

In this work we have built upon the ADMM as a fast and flexible solver for constrained optimization problems.
We have extended it in two directions. 
First, we allow for multiple constraint functions to be applied.
Unlike the previously proposed SDMM approach by \citet{Combettes2009}, we do not require that linear operators of the constraint functions, which may be needed for an efficient proximal operator formulation, be invertible.
Second, we address the case of a function that is convex in multiple arguments through an inexact block optimization method.
The proposed method, bSDMM, is effective in a range of constrained optimization problems that cannot fully be solved with e.g. proximal gradient methods.
As a result of its ADMM heritage, it is particularly suitable for applications where a fast and approximate solution is more important than an accurate one.

We showed its effectiveness to solve a hyperspectral unmixing problem, under the assumption of the linear mixing model, by performing a Non-Negative Matrix Factorization with multiple non-trivial constraints.
In a future work (Melchior et al., in prep.) we will utilize the bSDMM-NMF approach to separate stars and galaxies in astronomical images, a similar additive mixing problem, which requires several restrictive constraints for adequate performance.

As we believe in the usefulness of the algorithm and want to endorse reproducible research, we release the {\sf python} implementation of the algorithms presented here and the bSDMM-NMF as an open-source package at \url{https://github.com/pmelchior/proxmin}.

\section*{Acknowledgements}
We would like to thank Robert Vanderbei and Jonathan Eckstein for useful discussions regarding the algorithm, and Jim Bosch and Robert Lupton for comments on its astrophysical applications.

\bibliographystyle{spbasic}
\bibliography{references.bib}

\end{document}